\documentclass[article]{amsart}
\usepackage {amsmath, amssymb, color, textcomp}
\usepackage{mathtools}
\DeclarePairedDelimiter{\floor}{\lfloor}{\rfloor}

\newlength{\extramargin}

\makeatletter

\newtheorem{thm}{Theorem}[section]
\newtheorem{cor}[thm]{Corollary}
\newtheorem{lem}[thm]{Lemma}
\newtheorem{defn}[thm]{Definition}
\newtheorem{prop}[thm]{Proposition}

\newcommand{\eqnref}[1]{~(\ref{#1})}

\newtheorem{preremark}[thm]{Remark}
\newenvironment{remark}%
  {\begin{preremark}\upshape}{\end{preremark}}
\newtheorem{preexample}[thm]{Example}
  {\begin{preexample}\upshape}{\end{preexample}}

\numberwithin{equation}{section}

\numberwithin{equation}{section}

\numberwithin{thm}{section}
\newtheorem{lem/defn}[thm]{Lemma/Definition}
\newtheorem{preex/defn}[thm]{Example/Definition}
\newenvironment{ex/defn}%
  {\begin{preex/defn}\upshape}{\end{preex/defn}}

\newcommand{\ten}{\otimes}

\newcommand{\abs}[1]{\lvert#1\rvert}

\DeclareMathOperator{\End}{End}

\begin{document}

\title[Virasoro representations with central charges  $\frac{1}{2}$ and 1 on the Fock space $\mathit{F^{\otimes \frac{1}{2}}}$]{Virasoro representations with central charges  $\frac{1}{2}$ and 1 on the real neutral fermion Fock space $\mathit{F^{\otimes \frac{1}{2}}}$}

\author{Iana I. Anguelova}

\address{Department of Mathematics,  College of Charleston,
Charleston SC 29424 }
\email{anguelovai@cofc.edu}


\begin{abstract}
We study a family of  fermionic oscillator  representations of the Virasoro algebra via    2-point-local Virasoro fields  on the Fock space $\mathit{F^{\otimes \frac{1}{2}}}$ of a  neutral (real) fermion. We obtain the decomposition of  $\mathit{F^{\otimes \frac{1}{2}}}$ as a direct sum of  irreducible  highest weight Virasoro modules with central charge $c=1$.  As a corollary we obtain the decomposition of  the irreducible highest weight Virasoro modules with central charge $c=\frac{1}{2}$ into irreducible highest weight Virasoro modules with central charge $c=1$. As an application we show how  positive sum (fermionic) character formulas for the Virasoro modules of charge $c=\frac{1}{2}$  follow from these decompositions.
\end{abstract}

\maketitle

\section{Introduction}
\label{sec:intro}
This paper is a continuation of  \cite{Ang-BFVir}, and is a part of  a series   studying various particle correspondences in 2 dimensional conformal field theory from the point of view of chiral algebras (vertex algebras) and representation theory.
 The study of fermionic oscillator representations of the Virasoro algebra ($Vir$) is long-standing, dating back to \cite{Frenkel-BF}, \cite{FeiginFuks}, \cite{GoddardOliveSch}, \cite{GoddardOliveNahm}, \cite{GoddardKentOlive}, \cite{KacRaina}, \cite{Triality}, and many others. What was new in \cite{Ang-BFVir} is that we used 2-point-local  (local at both $z=w$ and $z=-w$ as opposed to just local at $z=w$) Virasoro fields generating the fermionic oscillator representations. In particular, in \cite{Ang-BFVir} we constructed  a 2-parameter family (depending on parameters $\lambda, b\in \mathbb{C}$) of   2-point-local   Virasoro fields with central charge $-2+12\lambda -12\lambda^2$  on the fermionic Fock space $\mathit{F^{\otimes \frac{1}{2}}}$ of the real fermion. In this paper we study the nature of these representations depending on the parameters,  and  their decomposition into irreducible modules. We show that for particular choices of the parameters $(\lambda, b)$  these two-point local Virasoro field representations   can produce each one of the well known discrete series of Virasoro representations.  Another  important and interesting particular case of the parameters is  that of $\lambda=\frac{1}{2}$, i.e., central charge $1$, and  $b\in \frac{1}{\sqrt{2}}\mathbb{Z}$. For this choice of $(\lambda, b)$ we obtain the decomposition of $\mathit{F^{\otimes \frac{1}{2}}}$ into irreducible highest weight modules for $Vir$ of central charge $1$. It is well known (going back to D. Friedan and I. Frenkel) that $\mathit{F^{\otimes \frac{1}{2}}}$ is a natural  highest weight module for $Vir$ of central charge $\frac{1}{2}$ via a one-point local Virasoro field, and as such it decomposes into two irreducible highest weight modules for $Vir$ of central charge $\frac{1}{2}$ (these are well known as two of the Ising minimal models). We show how each of these irreducible $Vir$ highest weight modules  of central charge $\frac{1}{2}$ decomposes into irreducible $Vir$  highest weight  modules  of central charge $1$. As an application, this allows us to directly write a positive sum (fermionic) character formula for these irreducible highest weight modules for $Vir$ of central charge $\frac{1}{2}$.
 There has been extensive research on  positive sum (fermionic) Virasoro character formulas, see for instance   \cite{Rinat1}, \cite{Foda1}, \cite{Bytsko}, \cite{Welsh},  \cite{FeiginFoda}. Fermionic-type character formulas are known for all the Virasoro minimal models (the discrete Virasoro series), although the vertex algebra (field theory) foundation for such character formulas is still lacking in the general case.  It is new here that we obtain these specific fermionic  formulas as a direct result of the decomposition of the charge $\frac{1}{2}$ "1-point-local" modules into charge 1 "2-point-local" modules (i.e., we explicitly use multi-local fields and twisted vertex algebra techniques). We hope that such approach is possible also more generally for the discrete Virasoro series by using quasi-particle Fock spaces in the place  of a single real fermionic Fock space of central charge  $\frac{1}{2}$. After bosonising such quasi-particle Fock spaces, the idea is then  to  obtain a decomposition of  the irreducible modules representing the one-point local Virasoro field of discrete-series-type   central charge into irreducible modules representing multi-local Virasoro fields of charge $1$.

\section{Notation and background}
\label{section:background}

We use the term "field" to mean  a series of the form
\begin{equation}
a(z)=\sum_{n\in \mathbf{Z}}a_{(n)}z^{-n-1}, \ \ \ a_{(n)}\in
\End(V), \ \ \text{such that }\ a_{(n_v)}v=0 \ \ \text{for any}\ v\in V, \ n_v\gg 0.
\end{equation}
Let
\begin{equation}
a(z)_-:=\sum_{n\geq 0}a_nz^{-n-1},\quad a(z)_+:=\sum_{n<0}a_nz^{-n-1}.
\end{equation}
\begin{defn} \label{defn:normalorder}{\bf (Normal ordered product)}
Let $a(z), b(z)$ be fields on a vector space $V$. Define
\begin{equation}
:\!a(z)b(w)\!:=a(z)_+b(w)+(-1)^{p(a)p(b)}b(w)a_-(z).
\end{equation}
One calls this the normal ordered product of $a(z)$ and $b(w)$.
\end{defn}
\begin{remark}
Let  $a(z), b(z)$ be fields on a vector space $V$. Then
$:\!a(z)b(\lambda z)\!:$ and $:\!a(\lambda z)b( z)\!:$ are well defined fields on $V$ for any $\lambda \in \mathbb{C}^*$.
\end{remark}
\begin{defn}(\cite{ACJ})  \label{defn:parity} {\bf ($2$-point-local fields) }
We say that a field $a(z)$ on a vector space $V$ is {\bf even} and $2$-point self-local at $(1; -1)$,  if there exist $n_0, n_1\in \mathbb{N}$ such that
\begin{equation}
(z- w)^{n_{0}}(z+w)^{n_{1}}[a(z),a(w)] =0.
\end{equation}
In this case we set the {\bf parity} $p(a(z))$ of $a(z)$ to be $0$.
\\
We set $\{a, b\}  =ab +ba$.We say that a field $a(z)$ on $V$ is $2$-point self-local at $(1; -1)$
and {\bf odd} if there exist $n_0, n_1\in \mathbb{N}$ such that
\begin{equation}
(z- w)^{n_{0}}(z+ w)^{n_{1}}\{a(z),a(w)\}=0.
\end{equation}
In this case we set the {\bf parity} $p(a(z))$ to be $1$. For brevity we will just write $p(a)$ instead of $p(a(z))$.\\
Finally,  if $a(z), b(z)$ are fields on $V$, we say that $a(z)$ and $b(z)$ are {\it $2$-point mutually local} at  $(1; -1)$
if there exist $n_0, n_1 \in \mathbb{N}$ such that
\begin{equation}
(z- w)^{n_{0}}(z+ w)^{n_{1}}\left(a(z)b(w)-(-1)^{p(a)p(b)}b(w)a(z)\right)=0.
\end{equation}
\end{defn}
For a rational function $f(z,w)$,  with poles only at $z=0$,  $z=\pm  w$, we denote by $i_{z,w}f(z,w)$
the expansion of $f(z,w)$ in the region $\abs{z}\gg \abs{w}$ (the region in the complex $z$ plane outside  the points $z=\pm w$), and correspondingly for
$i_{w,z}f(z,w)$.
The mathematical background of the well-known and often used (both in physics and in mathematics) notion of Operator Product Expansion (OPE) of product of two fields for case of usual locality ($N=1$) has been established for example in \cite{Kac}, \cite{LiLep}.
The following lemma extended the mathematical background  to the case of  $2$-point locality (in fact to $N$-point locality, for $N\in \mathbb{N}$):
\begin{lem} (\cite{ACJ}) {\bf (Operator Product Expansion (OPE))}\label{lem:OPE} \\
 Let $a(z)$, $b(w)$ be {\it $2$-point mutually local}. Then exists fields $c_{jk}(w)$, $j=0, 1; k=0, \dots , n_j-1$, such that we have
 \begin{equation}
 \label{eqn:OPEpolcor}
 a(z)b(w) =i_{z, w} \sum_{k=0}^{n_0-1}\frac{c_{0k}(w)}{(z-w)^{k+1}} + i_{z, w}\sum_{k=0}^{n_1-1}\frac{c_{1k}(w)}{(z+w)^{k+1}} +:\!a(z)b(w)\!:.
 \end{equation}
We call the fields $c_{jk}(w)$, $j=0, 1; k=0, \dots , n_j-1$ OPE coefficients. We will write the above OPE as
 \begin{equation}
 a(z)b(w) \sim  \sum_{k=0}^{n_0-1}\frac{c_{0k}(w)}{(z-w)^{k+1}} + \sum_{k=0}^{n_1-1}\frac{c_{1k}(w)}{(z+w)^{k+1}}.
 \end{equation}
 \end{lem}
  The $\sim $ signifies that we have only written the singular part, and also we have omitted writing explicitly the expansion $i_{z, w}$, which we do acknowledge  tacitly. Often also the following notation is used for short:
 \begin{equation}\label{contraction}
\lfloor
ab\rfloor=a(z)b(w)-:a(z)b(w):= [a(z)_-,b(w)],
\end{equation}
i.e.,  the {\it contraction} of any two fields
$a(z)$ and $b(w)$ is in fact also the $i_{z, w}$ expansion of the singular part of the OPE of the two fields $a(z)$ and $b(w)$.

The OPE expansion in the multi-local case allowed us to  extend the Wick's Theorem (see e.g., \cite{MR85g:81096}, \cite{MR99m:81001}) to the case of multi-locality (see \cite{ACJ}), and we will use it in what follows, together with the Taylor Expansion Lemma (see \cite{ACJ}).

\section{The Fock space  $\mathit{F^{\ten \frac{1}{2}}}$ and 2-point-local field representations of  $Vir$}
\label{section:main}

 We recall the definitions and notations for the Fock space $\mathit{F^{\ten \frac{1}{2}}}$ as in \cite{Frenkel-BF}, \cite{DJKM3}, \cite{Kac-Lie}, \cite{WangDuality}; in particular we follow the notation of \cite{WangDuality}, \cite{WangDual}.

Consider a single odd self-local field $\phi ^D(z)$, which we index in the form $\phi ^D(z)=\sum _{n\in \mathbb{Z}+\frac{1}{2}} \phi^D_n z^{-n-\frac{1}{2}}$.
The OPE of $\phi ^D(z)$ is given by
\begin{equation}
\label{equation:OPE-D}
\phi ^D(z)\phi ^D(w)\sim \frac{1}{z-w}.
\end{equation}
This OPE completely determines the commutation relations between the modes $\phi^D_n$, $n\in \mathbb{Z} +\frac{1}{2}$:
\begin{equation}
\label{eqn:Com-D}
\{\phi^D_m,\phi^D_n\}:=\phi^D_m\phi^D_n + \phi^D_n\phi^D_m=\delta _{m, -n}1.
\end{equation}
and so the modes generate a Clifford algebra $\mathit{Cl_D}$. The field $\phi ^D(z)$ is usually  called a ``real neutral fermion field".
The Fock space, denoted  by $\mathit{F^{\ten \frac{1}{2}}}$,  of the  field $\phi ^D(z)$ is the highest weight module of $\mathit{Cl_D}$ with vacuum vector $|0\rangle $, so that $\phi^D_n|0\rangle=0 \ \text{for} \  n >0$.   This well known Fock space is often called the Fock space of the free real neutral fermion (see e.g.   \cite{DJKM6}, \cite{Triality}, \cite{Wang}, \cite{WangKac},  \cite{WangDual},  \cite{WangDuality}, \cite{Rehren}).
$\mathit{F^{\ten \frac{1}{2}}}$  has basis
\begin{equation}
\{ \phi^D_{-n_k-\frac{1}{2}}\dots \phi^D_{-n_2-\frac{1}{2}}\phi^D_{-n_1-\frac{1}{2}}|0\rangle, |0\rangle |\ \text{where}\  n_k>\dots >n_2>n_1\geq 0, n_i\in \mathbb{Z}_{\geq 0}, i=1, 2, \dots, k\}
\end{equation}
We recall the various gradings on $\mathit{F^{\ten \frac{1}{2}}}$.
The space $\mathit{F^{\ten \frac{1}{2}}}$ has  a well known $\mathbb{Z}_2$ grading given by  $ k\  mod\  2$,
\[
\mathit{F^{\ten \frac{1}{2}}}=\mathit{F_{\bar{0}}^{\ten \frac{1}{2}}}\oplus \mathit{F_{\bar{1}}^{\ten \frac{1}{2}}},
\]
where $\mathit{F_{\bar{0}}^{\ten \frac{1}{2}}}$ (resp. $\mathit{F_{\bar{1}}^{\ten \frac{1}{2}}}$) denote the even (resp. odd) components of $\mathit{F^{\ten \frac{1}{2}}}$. This $\mathbb{Z}_2$ grading can be extended  to a $\mathbb{Z}_{\geq 0}$ grading $\tilde{L}$, called ``length", by setting
\begin{equation}
\tilde{L} (\phi^D_{-n_k-\frac{1}{2}}\dots \phi^D_{-n_2-\frac{1}{2}}\phi^D_{-n_1-\frac{1}{2}}|0\rangle)=k.
\end{equation}
Using the $\mathbb{Z}_2$ grading the space $\mathit{F^{\ten \frac{1}{2}}}$ can be given a super-vertex algebra structure, as is known from e.g. \cite{Triality}, \cite{Wang}, \cite{Kac}.

In \cite{ACJ2} and \cite{Ang-BFVir} we introduced a $\mathbb{Z}$ grading $dg$ on $\mathit{F^{\ten \frac{1}{2}}}$ by assigning $dg(|0\rangle)=0$ and
\begin{equation*}
dg(\phi^D_{-n_k-\frac{1}{2}}\dots \phi^D_{-n_2-\frac{1}{2}}\phi^D_{-n_1-\frac{1}{2}}|0\rangle)=\#\{i=1, 2, \dots, k|n_i=\text{odd}\}\ -\#\{i=1, 2, \dots, k|n_i=\text{even}\}.
\end{equation*}
Denote the homogeneous component of degree $dg=n \in \mathbb Z$ by $\mathit{F_{(n)}^{\ten \frac{1}{2}}}$, hence as vector spaces we have
\begin{equation}
\mathit{F^{\ten \frac{1}{2}}} = \oplus _{n\in \mathbb{Z}} \mathit{F_{(n)}^{\ten \frac{1}{2}}}.
\end{equation}
We define the special vectors $v_n\in  \mathit{F_{(n)}^{\ten \frac{1}{2}}}$  by
\begin{align}
&v_0=|0\rangle \in  \mathit{F_{(0)}^{\ten \frac{1}{2}}};\\
&v_n=\phi^D_{-2n+1-\frac{1}{2}}\dots \phi^D_{-3-\frac{1}{2}}\phi^D_{-1-\frac{1}{2}}|0\rangle \in  \mathit{F_{(n)}^{\ten \frac{1}{2}}}, \quad \text{for}\ n>0;\\
&v_{-n}=\phi^D_{-2n+2-\frac{1}{2}}\dots \phi^D_{-2-\frac{1}{2}}\phi^D_{-\frac{1}{2}}|0\rangle\in  \mathit{F_{(-n)}^{\ten \frac{1}{2}}}, \quad \ n>0.
\end{align}
Note that the vectors $v_n\in \mathit{F_{(n)}^{\ten \frac{1}{2}}}$ have minimal length $\tilde{L}=|n|$  among the vectors within $\mathit{F_{(n)}^{\ten \frac{1}{2}}}$, and they are in fact the unique (up-to a scalar) vectors minimizing  the length $\tilde{L}$, such that the index $n_k$ is minimal too ($n_k$ is identified from the smallest  index appearing). One can view each of the vectors $v_n$ as a vacuum-like vector in $\mathit{F_{(n)}^{\ten \frac{1}{2}}}$, see below,  and the $dg$ grading as the analogue of the "charge" grading.

Note also that the super-grading derived from the second grading $dg$ is compatible with the super-grading derived from the first grading $\tilde{L}$, as  $dg \ mod\  2 =\tilde{L}\  mod\  2$.

Lastly, we  recall the grading $deg_h$ (from \cite{Ang-BFVir})) on each of the components $\mathit{F_{(n)}^{\ten \frac{1}{2}}}$ for each $n\in \mathbb{Z}$. Consider a monomial vector $v=\phi^D_{-n_k-\frac{1}{2}}\dots \phi^D_{-n_2-\frac{1}{2}}\phi^D_{-n_1-\frac{1}{2}}|0\rangle$ from $\mathit{F_{(n)}^{\ten \frac{1}{2}}}$. One can view this vector as an "excitation" from the vacuum-like vector $v_n$, and count the $n_i$ that  should have been in $v$ as compared to $v_n$, and also the  $n_i$ that  should not have been in $v$ as compared to $v_n$. Thus the grading $deg_h$ (one can think of it as "energy") is defined as
\begin{align*}
deg_h (v)& =\sum \{\floor[\big ]{\frac{n_l+1}{2}}\  |\  n_l\  \text{that \ should \ have \ been \ there\ but \ are \ not}\}\\
&\quad +\sum \{\floor[\big ]{\frac{n_l+1}{2}}\ |\ n_l\ \text{that \ should \ not have \ been \ there\ but\ are}\};
\end{align*}
here $\floor[\big ]{x}$ denotes the floor function, i.e.,   $\floor[\big ]{x}$ denotes the largest integer smaller or equal to  $x$.  Denote by  $\mathit{F_{(n, k)}^{\ten \frac{1}{2}}}$ the linear span of all vectors of grade $deg_h =k$ in $\mathit{F_{(n)}^{\ten \frac{1}{2}}}$.

Recall the Heisenberg algebra $\mathcal{H}_{\mathbb{Z}}$ with relations $[h_m,h_n]=m\delta _{m+n,0}1$, \ $m,n$   integers.
We proved in \cite{Ang-BFVir} (and by the less traditional bicharacter construction in \cite{AngTVA}) that each $\mathit{F_{(n)}^{\ten \frac{1}{2}}}$ is an irreducible highest weight module for the Heisenberg algebra via a 2-point-local field:
\begin{prop} \label{prop:heisdecomp} (\cite{AngTVA}, \cite{Ang-BFVir})
The 2-point-local field  $h^D(z)$ given by:
\begin{equation}
\label{eqn:normal-order-h-D}
h^D (z)= \frac{1}{2}:\phi ^D(z)\phi ^D(-z) :
\end{equation}
has only  odd-indexed modes ($h^D (z)=-h^D (-z)$),   $h^D (z)=\sum _{n\in \mathbb{Z}} h_{n} z^{-2n-1}$,  and has OPE with itself given by:
\begin{equation}
\label{eqn:HeisOPEsD}
h^D (z)h^D (w)\sim \frac{zw}{(z^2-w^2)^2} \sim \frac{1}{4}\frac{1}{(z-w)^2} - \frac{1}{4}\frac{1}{(z+w)^2}.
\end{equation}
 Hence its  modes, $h_n, \ n\in \mathbb{Z}$, generate the Heisenberg algebra $\mathcal{H}_{\mathbb{Z}}$.
The neutral real fermion Fock space  $\mathit{F^{\ten \frac{1}{2}}}$ is thus a module for $\mathcal{H}_{\mathbb{Z}}$ via this 2-point-local field representation and decomposes into irreducible highest weight modules for  $\mathcal{H}_{\mathbb{Z}}$ under this action  as follows:
\begin{equation}
\mathit{F^{\ten \frac{1}{2}}} = \oplus _{m\in \mathbb{Z}} \mathit{F_{(m)}^{\ten \frac{1}{2}}}\cong \oplus _{m\in \mathbb{Z}} B_m,  \quad \text{where} \quad \mathit{F_{(m)}^{\ten \frac{1}{2}}}\cong B_m,\quad
B_m \cong  \mathbb{C}[x_1, x_2, \dots , x_n, \dots ],\quad \forall \ m\in \mathbb{Z}.
\end{equation}
\end{prop}

Recall the  well-known  Virasoro algebra $Vir$, the central extension of the complex polynomial vector fields on the circle. The Virasoro  algebra $Vir$ is the Lie algebra with generators $L_n$, $n\in \mathbb{Z}$, and central element $C$, with
commutation relations
\begin{equation}
\label{eqn:VirCRs}
[L_m, L_n] =(m-n)L_{m+n} +\delta_{m, -n}\frac{(m^3-m)}{12}C; \quad [C, L_m]=0, \ m, n\in \mathbb{Z}.
\end{equation}
Equivalently, the 1-point-local Virasoro field
$L(z): =\sum _{n\in \mathbb{Z}} L_{n} z^{-n-2}$
has OPE with itself given by:
\begin{equation}
\label{eqn:VirOPEs}
L(z)L(w)\sim \frac{C/2}{(z-w)^4} + \frac{2L(w)}{(z-w)^2}+ \frac{\partial_{w}L(w)}{(z-w)}.
\end{equation}

We denote by $M(c,h)$ the irreducible highest weight  $Vir$ module  with central charge $c\in \mathbb{C}$, where $h\in \mathbb{C}$ is the weight of the operator $L_0$ acting on the highest weight vector $v_h$, i.e., $M(c,h)$ is an irreducible  $Vir$ module  generated by a vector $v_h\in M(c,h)$ with $Cv_h=cv_h$, $L_n v_h=0$  for any $n> 0$, and  $L_0 v_h =hv_h$. (Note that in some of the literature such  $Vir$ modules are called "lowest weight", instead of "highest weight", here we follow the convention of Kac, as in \cite{Kacdet}, \cite{KacRaina}, \cite{Kac-Lie}). The seminal paper \cite{FeiginFuks} (using the Kac determinant formula, \cite{Kacdet})  delineates the cases where $M(c,h)$ is a quotient by a nontrivial factor of the corresponding Verma module $VM(c,h)$, as opposed to the "generic" cases where $M(c,h)\cong VM(c,h)$.

Denote the formal character of a highest weight module $V$ with highest weight $(c, h)$ to be
\[
ch_q^{Vir}V:= tr_V q^{L_0} :=\sum _{j\in \mathbb{Z}_+}(dim V_{h+j})q^{h+j},
\]
where $V_{h+j}$ is the eigenspace of $L_0$ of weight $h+j$.

It is well known that  $\mathit{F^{\ten \frac{1}{2}}}$  is  a module for the Virasoro algebra with central charge $c=\frac{1}{2}$ (see for example \cite{Triality}, \cite{Wang}, \cite{WangDual}) via the   1-point local Virasoro field $L^{1/2}(z)$  given by
\begin{equation}
\label{eqn:Vir1/2}
L^{1/2}(z)=\frac{1}{2}: \partial_z \phi ^D (z) \phi ^D(z):.
\end{equation}
Furthermore, it is well known that as Vir modules
\[
\mathit{F_{\bar{0}}^{\ten \frac{1}{2}}}\cong M\left(\frac{1}{2}, 0\right); \quad \mathit{F_{\bar{1}}^{\ten \frac{1}{2}}}\cong M\left(\frac{1}{2}, \frac{1}{2}\right); \quad \mathit{F^{\ten \frac{1}{2}}}\cong M\left(\frac{1}{2}, 0\right)\oplus M\left(\frac{1}{2}, \frac{1}{2}\right).
\]
In \cite{Ang-BFVir} we proved the following Jacobi identity holds:
\begin{prop}\label{prop:Jacobi}(Corollary to Proposition \ref{prop:heisdecomp})
Define the graded dimension (character) of the Fock space  $\mathit{F}=\mathit{F^{\ten \frac{1}{2}}}$ as
\[
ch \mathit{F} :=tr_{\mathit{F}}q^{L^{1/2}_0}z^{h_0}.
\]
We have
\begin{align}
\label{eqn:jac1}
ch \mathit{F}&=\prod_{i=1}^{\infty} (1+zq^{2i-1+\frac{1}{2}})(1+z^{-1}q^{2i-2+\frac{1}{2}})\\
\label{eqn:jac2}
&=\frac{1}{\prod_{i=1}^{\infty} (1-q^{2i})}\sum_{n\in \mathbb{Z}} z^nq^{\frac{n}{2}}q^{n^2}.
\end{align}
By comparing the two identities we get the Jacobi identity
\begin{equation}
\label{eqn:Jacid}
\prod_{i=1}^{\infty} (1-q^{2i})(1+zq^{2i-\frac{1}{2}})(1+z^{-1}q^{2i-\frac{3}{2}})=\sum_{m\in \mathbb{Z}}z^m q^{\frac{m(2m+1)}{2}}.
\end{equation}
\end{prop}
Observe that we can specialize the graded dimension $ch \mathit{F}$ to the Virasoro character by evaluating at $z=1$:
\begin{cor}
\begin{align}
\label{eqn:jac1Vir}
ch_q^{Vir} \mathit{F}=ch_q M\left(\frac{1}{2}, 0\right) + ch_q M\left(\frac{1}{2}, \frac{1}{2}\right)&=\prod_{i=0}^{\infty} (1+q^{\frac{1}{2}}q^{i})\\
\label{eqn:jac2Vir}
&=\frac{1}{\prod_{i=1}^{\infty} (1-q^{2i})}\sum_{n\in \mathbb{Z}} q^{n^2 +\frac{n}{2}}.
\end{align}
\end{cor}
If we use the $q$-Pochhammer symbol notation,
\begin{equation*}
(a; q)_{\infty}: =\prod_{i=0}^{\infty} (1-aq^{i}); \quad (a; q)_m: =\prod_{i=0}^{m-1} (1-aq^{i}), \quad \text{with} \ (a; q)_0: =1;
\end{equation*}
we can rewrite the first formula \eqref{eqn:jac1Vir}
\begin{equation*}
ch_q^{Vir} \mathit{F} = ch_q M\left(\frac{1}{2}, 0\right) + ch_q M\left(\frac{1}{2}, \frac{1}{2}\right)= (-q^{\frac{1}{2}}; q)_{\infty}.
\end{equation*}
The  formula  \eqref{eqn:jac1Vir} is well known (going back to   I. Frenkel and D. Friedan, as well as \cite{FeiginFuks}), but here we obtain it as a direct evaluation at $z=1$ of the formula \eqref{eqn:jac1}.
We can refine the proof of \eqref{eqn:jac2} from \cite{Ang-BFVir} to obtain separate formulas for $ch_q M\left(\frac{1}{2}, 0\right)$ and $ch_q M\left(\frac{1}{2}, \frac{1}{2}\right)$:
\begin{prop}\label{Virchar}
\begin{align}
\label{eqn:char0Vir}
ch_q M\left(\frac{1}{2}, 0\right)&=\frac{1}{\prod_{i=1}^{\infty} (1-q^{2i})}\sum_{\substack{n\in \mathbb{Z}\\ n\ \text{even}}} q^{n^2 +\frac{n}{2}};\\
\label{eqn:char1Vir}
ch_q M\left(\frac{1}{2}, \frac{1}{2}\right)&=\frac{1}{\prod_{i=1}^{\infty} (1-q^{2i})}\sum_{\substack{n\in \mathbb{Z}\\ n\ \text{odd}}} q^{n^2 +\frac{n}{2}}
\end{align}
\end{prop}
These formulas have all positive coefficients in their sum, as
\[
\frac{1}{\prod_{i=1}^{\infty} (1-q^{i})}=\sum_{k\geq 0} p(k)q^k,
 \]
 where $p(k)$  is  the number of  partitions $0<k_1\leq \dots \leq k_l$ of $k=k_1+\dots +k_l$. Such formulas as \eqref{eqn:char0Vir} and \eqnref{eqn:char1Vir} are called fermionic (see e.g.  \cite{Rinat1}), to distinguish them from the character sums with alternating sign coefficients (as in \cite{RochaCharidi}), which are called bosonic.\\
\emph{Proof:} We have, as vector spaces, that
\begin{equation*}
M\left(\frac{1}{2}, 0\right)\cong \mathit{F_{\bar{0}}^{\ten \frac{1}{2}}} =  \oplus_{\substack{n\in \mathbb{Z}\\ n\ \text{even}}} \mathit{F_{(n)}^{\ten \frac{1}{2}}};\quad \quad
M\left(\frac{1}{2}, \frac{1}{2}\right)\cong \mathit{F_{\bar{1}}^{\ten \frac{1}{2}}}=\oplus_{\substack{n\in \mathbb{Z}\\ n\ \text{odd}}} \mathit{F_{(n)}^{\ten \frac{1}{2}}}.
\end{equation*}
We can now use the  decomposition  of each $\mathit{F_{(n)}^{\ten \frac{1}{2}}}$ into $\mathit{F_{(n, k)}^{\ten \frac{1}{2}}}$ by $deg_h$. From Proposition \ref{prop:heisdecomp} it follows that a  basis for $\mathit{F_{(n, k)}^{\ten \frac{1}{2}}}$ is given by the elements $v_{h, n, \vec{k}}=h_{-k_{l}}\dots h_{-k_{1}}v_n$ with indecies varying with partitions $0<k_1\leq \dots \leq k_l$ of $k=k_1+\dots +k_l$.
First, we have
\begin{align*}
L^{1/2}_0 v_0& = L^{1/2}_0 |0\rangle =0\\
L^{1/2}_0 v_n &= L^{1/2}_0 \phi^D_{-2n+1-\frac{1}{2}}\dots \phi^D_{-3-\frac{1}{2}}\phi^D_{-1-\frac{1}{2}}|0\rangle \\
&\quad = \left((1+\frac{1}{2})+ (3+\frac{1}{2}) +\dots +(2n-1+\frac{1}{2})\right) v_n=(n^2+\frac{n}{2}) v_n, \quad \text{for}\ n>0;\\
L^{1/2}_0 v_{-n} &= L^{1/2}_0 \phi^D_{-2n+2-\frac{1}{2}}\dots \phi^D_{-2-\frac{1}{2}}\phi^D_{-\frac{1}{2}}|0\rangle \\
&\quad =\left((0+\frac{1}{2}) + (2+\frac{1}{2}) + (4+\frac{1}{2})+\dots +(2n-2+\frac{1}{2})\right) v_{-n} =(n^2-\frac{n}{2})v_{-n}, \quad \text{for}\ n>0.
\end{align*}
Now a direct calculation shows that (\cite{Ang-BFVir}):
\begin{equation}
\label{eqn:h-normorder}
:h^D(w)^2:= :h^D(w)h^D(w):=\frac{1}{4}:(\partial_{-w}\phi ^D(-w))\phi ^D(-w): + \frac{1}{4}:(\partial_{w}\phi ^D(w))\phi ^D(w): -\frac{1}{2w}h^D(w).
\end{equation}
We can calculate by direct use of Wick's Theorem the OPE between  $:h^D(z)^2:$ and $h^D(w)$, and thus by use of the equation above the commutator of  $L^{1/2}_0$ and $h^D(w)$; and we obtain that
\begin{equation}
[L^{1/2}_0, h^D_k]=-2k h^D_k.
\end{equation}
Hence
\begin{align*}
L^{1/2}_0 v_{h, 0, \vec{k}}& =L^{1/2}_0 h_{-k_{l}}\dots h_{-k_{1}}v_n=\left(2k_1+\dots +2k_l\right)v_{h, 0, \vec{k}};\\
L^{1/2}_0 v_{h, n, \vec{k}}& =L^{1/2}_0 h_{-k_{l}}\dots h_{-k_{1}}v_n=\left(2k_1+\dots +2k_l +n^2+\frac{n}{2}\right)v_{h, n, \vec{k}} \quad \text{for}\ n>0;\\
L^{1/2}_0 v_{h, -n, \vec{k}}& =L^{1/2}_0 h_{-k_{l}}\dots h_{-k_{1}}v_{-n}=\left(2k_1+\dots +2k_l +n^2 -\frac{n}{2}\right)v_{h, -n, \vec{k}} \quad \text{for}\ n>0.
\end{align*}
Now since there are $p(k)$ such elements for partitions $0<k_1\leq \dots \leq k_l$ of $k=k_1+\dots +k_l$, we have for $\mathit{F_{\bar{0}}^{\ten \frac{1}{2}}}\cong M\big(\frac{1}{2}, 0\big)$
\begin{align*}
ch_q M\big(\frac{1}{2}, 0\big)&= tr_{\mathit{F_{\bar{0}}^{\ten \frac{1}{2}}}} q^{L^{1/2}_0}=
\sum_{k\geq 0}p(k)q^{2k}+\sum_{\substack{n, k\in \mathbb{Z}_+\\n\ even}}p(k)\left(q^{2k+n^2+\frac{n}{2}} +q^{2k+n^2-\frac{n}{2}}\right)\\
& =\sum_{k\in \mathbb{Z}, k\geq 0}p(k)q^{2k}\cdot \left(1+\sum_{\substack{n\in \mathbb{Z}_+\\n\ even}} \left(q^{n^2+\frac{n}{2}} +q^{n^2-\frac{n}{2}}\right)\right)\\
&=\frac{1}{\prod_{i=1}^{\infty} (1-q^{2i})}\sum_{\substack{n\in \mathbb{Z}\\n\ even}} q^{n^2+\frac{n}{2}}.
\end{align*}
The calculation for  $\mathit{F_{\bar{1}}^{\ten \frac{1}{2}}}\cong M\big(\frac{1}{2}, \frac{1}{2}\big)$ is even simpler. $\hspace{7cm} \square$

We now turn to the new representations in this paper, the 2-point local Virasoro field representations on the Fock space  $\mathit{F^{\ten \frac{1}{2}}}$.
In \cite{Ang-BFVir} we proved that besides $L^{1/2}(z)$ there is a 2-parameter family of 2-point-local Virasoro fields on $\mathit{F^{\ten \frac{1}{2}}}$:
\begin{prop}\label{prop:VirasoroNew} (\cite{Ang-BFVir})
The 2-point-local field
\begin{equation}
\label{eqn:normal-order-L-D}
L^{1}(z^2): =\frac{1}{2z^2}:h^D (z)h^D (z):
\end{equation}
has only  even-indexed modes,    $L^{1}(z^2):=\sum _{n\in \mathbb{Z}} L^1_{n} (z^2)^{-n-2}$ and  its modes $L_n$ satisfy the Virasoro algebra commutation relations with central charge $c=1$:
\[
[L_m, L_n] =(m-n)L_{m+n} +\delta_{m, -n}\frac{(m^3-m)}{12}.
\]
Equivalently, the 2-point-local field $L^{1}(z^2)$ has OPE with itself given by:
\begin{equation}
\label{eqn:VirOPEsD}
L^{1}(z^2)L^{1}(w^2)\sim \frac{1/2}{(z^2-w^2)^4} + \frac{2L^{1}(w^2)}{(z^2-w^2)^2}+ \frac{\partial_{w^2}L^{1}(w^2)}{(z^2-w^2)}.
\end{equation}
Furthermore, the 2-point-local field
\begin{equation}
\label{eqn:normal-order-L-D-lambda}
L^{\lambda, b}(z^2): =L^{1}(z^2) + \frac{1 -2\lambda }{4z^2}\partial_z h^D (z)-\frac{b}{z^3}h^D(z) +\frac{\left(b+\frac{1 -2\lambda }{4}\right)^2 -4\left(\frac{1 -2\lambda }{4}\right)^2}{2z^4}=\sum _{n\in \mathbb{Z}} L^{\lambda, b}_{n} (z^2)^{-n-2}
\end{equation}
is a Virasoro field for every $\lambda, b \in \mathbb{C}$ with central charge $-12\lambda^2+12\lambda -2$. If $\lambda =\frac{1}{2}, \ b=0$, $L^{\frac{1}{2}, 0}(z^2)=L^{1}(z^2)$.
\end{prop}
We now turn to the types of $Vir$ representations these 2-point-local Virasoro fields generate on $\mathit{F^{\ten \frac{1}{2}}}$, depending on the choices of the parameters $(\lambda, b)$.
First, observe that  these representations of $Vir$,  although  of course expressible as   2-point local fermionic oscillator representations through the generating field $\phi ^D(z)$ (see \eqref{eqn:h-normorder}),  are in fact only dependent on the descendent  Heisenberg field $h^D(z)$ from Proposition \ref{prop:heisdecomp}.  It is then  immediate that each
$\mathit{F_{(n)}^{\ten \frac{1}{2}}}$ for $n\in \mathbb{Z}$ is also a submodule  for these  2-point-local $Vir$ representations, and is in  fact a highest weight module:
\begin{equation}
L^{\lambda, b}_m v_n=0, \quad \text{for any} \ m,n\in \mathbb{Z}, \quad \text{where} \ m>0,
\end{equation}
and we have  for the highest weight vectors $v_n$ ($n\in \mathbb{Z}$)
\begin{align*}
L^{\lambda, b}_0 v_n &=\left(\frac{n^2}{2}-bn -(\frac{1-2\lambda}{4})n+\frac{\left(b+\frac{1 -2\lambda }{4}\right)^2 -4(\frac{1 -2\lambda }{4})^2}{2}\right)v_n \\
&=\left(\frac{\left(b+\frac{1 -2\lambda}{4}-n\right)^2}{2} -2\left(\frac{1 -2\lambda }{4}\right)^2\right)v_n.
\end{align*}
Observe that
\[
c=1 \quad \text{only for\ }\ \lambda=\frac{1}{2}; \quad c=-12\lambda^2+12\lambda -2< 1 \quad \text{for\ }\ \lambda \in \mathbb{R}\setminus \{\frac{1}{2}\}.
\]
 For "generic" cases of real $0\leq c\leq 1$  and $h$ the Verma modules $VM(c,h)$  are  irreducible  representations (\cite{Kacdet}, \cite{FeiginFuks}), and for those "generic" $(c, h)$ we have  $\mathit{F_{(n)}^{\ten \frac{1}{2}}}\cong VM(c, h)=M(c, h)$.   Furthermore, one is of course interested in the cases where the representation is unitary, and that leaves only  the discrete series  (see \cite{FriedanQS}, and e.g. \cite{KacRaina}), or the case of  $c=1$ (see e.g. \cite{Kacdet}, \cite{FeiginFuks}, \cite{FriedanQS}, \cite{KacRaina}). For $\mathit{F_{(n)}^{\ten \frac{1}{2}}}$ to be   of    discrete series type  we need to have (\cite{Kacdet}, \cite{FeiginFuks}, \cite{FriedanQS}, \cite{KacRaina}):
\begin{align*}
c&=-12\lambda^2+12\lambda -2=1-\frac{6}{(m+2)(m+3)}, \quad \text{for}\ m\in \mathbb{Z}_{\geq 0}, \\
h&=h_{r, s}=\frac{\left((m+3)r-(m+2)s\right)^2-1}{4(m+2)(m+3)}  \quad \text{for}\ r, s\in \mathbb{Z}_+, \ 1\leq s\leq r\leq m+1;
\end{align*}
which here gives us
\begin{equation*}
\left(1 -2\lambda \right)^2=\frac{2}{(m+2)(m+3)}.
\end{equation*}
Now we have  to consider for such $\lambda$ also the highest weight $h$, which gives us
\begin{equation*}
\frac{\left(b+\frac{1 -2\lambda}{4}-n\right)^2}{2} =\frac{\left((m+3)r-(m+2)s\right)^2}{4(m+2)(m+3)}  \quad \text{for\ some}\ r, s\in \mathbb{Z}_+, \ 1\leq s\leq r\leq m+1;
\end{equation*}
and thus
\begin{equation*}
\left(\frac{b+\frac{1 -2\lambda}{4}-n}{\frac{1 -2\lambda}{4}}\right)^2=4\left((m+3)r-(m+2)s\right)^2, \quad \text{for}\ r, s\in \mathbb{Z}_+, \ 1\leq s\leq r\leq m+1.
\end{equation*}
Hence we can find $b\in \mathbb{R}$ which will produce a discrete-series-type representation on  $\mathit{F_{(n)}^{\ten \frac{1}{2}}}$. To summarize, if the parameters $(\lambda, b)$  satisfy
\begin{align*}
\left(2\lambda -1 \right)&=\pm \sqrt{\frac{2}{(m+2)(m+3)}}, \quad \text{for}\ m\in \mathbb{Z}_+;\\
\left(b-n\right) &=\pm\frac{2\left((m+3)r-(m+2)s\right)\pm 1}{2\sqrt{2(m+2)(m+3)}}, \quad \text{for}\ r, s\in \mathbb{Z}_+, \ 1\leq s\leq r\leq m+1;
\end{align*}
then the submodule $\mathit{F_{(n)}^{\ten \frac{1}{2}}}$  belongs to the discrete series. Thus we can choose values of the parameters $(\lambda, b)$ to produce a 2-point local fermionic oscillator  field representation of each of the discrete series Virasoro modules.

Next we  consider the  case of $c=1$, i.e.,  $\lambda =\frac{1}{2}$. In this case we either have
\begin{equation*}
L^{\frac{1}{2}, b}_0 v_n =\frac{(b-n)^2}{2}v_n=\frac{m^2}{4}v_n,    \quad \text{for\ some}\ m\in \mathbb{Z};
\end{equation*}
i.e.,
\begin{equation}
(b-n)^2=\frac{m^2}{2},    \quad \text{for\ some}\ m\in \mathbb{Z};
\end{equation}
or the submodule $\mathit{F_{(n)}^{\ten \frac{1}{2}}}$ is irreducible. It is clear that if for some $n_1\in \mathbb{Z}$ the parameter $b$ satisfies
\begin{equation}
(b-n_1)^2=\frac{m^2}{2},    \quad \text{for\ some}\ m\in \mathbb{Z};
\end{equation}
then for all other $n\in \mathbb{Z}$, $n\neq n_1$, we have  $(b-n)^2\neq \frac{m^2}{2}, \  \text{for\ any}\ m\in \mathbb{Z}$. Hence if the module $\mathit{F_{(n_1)}^{\ten \frac{1}{2}}}$ is indeed reducible, then for all other $n\in \mathbb{Z}$, $n\neq n_1$, $\mathit{F_{(n)}^{\ten \frac{1}{2}}}$ will be irreducible.  Thus, when $\lambda =\frac{1}{2}$ we either have $(b-n)^2\neq \frac{m^2}{2}, \  \text{for\ any}\ n, m\in \mathbb{Z}$ and thus all modules $\mathit{F_{(n)}^{\ten \frac{1}{2}}}$ are irreducible (we call such $b$ "generic"). Or exactly one the modules  $\mathit{F_{(n_1)}^{\ten \frac{1}{2}}}$ is completely reducible, the others are irreducible. The general structure of these reducible highest $Vir$ modules with central charge $1$ is well known from \cite{FeiginFuks} (see also \cite{KacRaina}), but here we can  actually describe explicitly the singular vectors generating the submodules.   We can  without loss of generality assume that for $\lambda =\frac{1}{2}$  the one reducible submodule is $\mathit{F_{(0)}^{\ten \frac{1}{2}}}$ (i.e., $b^2=\frac{m^2}{2}, \  \text{for\ some}\ m\in \mathbb{Z}$), and thus the other $\mathit{F_{(n)}^{\ten \frac{1}{2}}}$, $n\neq 0$,  are irreducible.

\begin{lem} Let $\lambda =\frac{1}{2}, \ b=\frac{m}{\sqrt{2}}, \  m\in \mathbb{Z}$ is fixed. \\
\textbf{Case I.} For $m\geq  0$, the following vectors in $\mathit{F_{(0)}^{\ten \frac{1}{2}}}$, indexed by  $k\in \mathbb{Z}$, $k\geq 0$,   exhaust all singular vectors for the two-point local $Vir$ representation on $\mathit{F_{(0)}^{\ten \frac{1}{2}}}$:
\begin{align*}
\tilde{v}_{m, 0}&= v_0=|0\rangle, \quad \text{for}\ k=0;\\
\tilde{v}_{m, k}&=\phi^D_{-2(k+m)+1-\frac{1}{2}}\phi^D_{-2(k+m-1)+1-\frac{1}{2}}\dots \phi^D_{-2m -1 -\frac{1}{2}}\phi^D_{-2(k-1)-\frac{1}{2}}\phi^D_{-2(k-2)-\frac{1}{2}}\dots \phi^D_{-\frac{1}{2}}|0\rangle, \quad \text{for}\ k>0.
\end{align*}
We have
\begin{align*}
L^{\frac{1}{2}, \frac{m^2}{2}}_j \tilde{v}_{m, k}&=0, \quad \text{for \ any} \ j>0, k\geq 0;\\
L^{\frac{1}{2}, \frac{m^2}{2}}_0 \tilde{v}_{m, k}&=\frac{1}{4}(m+2k)^2\tilde{v}_{0, k}.
\end{align*}
\textbf{Case II.} For $m < 0$, the following vectors  in $\mathit{F_{(0)}^{\ten \frac{1}{2}}}$, indexed by  $k\in \mathbb{Z}$, $k \geq -m$,  exhaust all singular vectors for the two-point local $Vir$ representation on $\mathit{F_{(0)}^{\ten \frac{1}{2}}}$:
\begin{align*}
\tilde{v}_{m, -m}&= v_0=|0\rangle,  \quad \text{for}\ k=-m;\\
\tilde{v}_{m, k}&=\phi^D_{-2(k+m)+1-\frac{1}{2}}\phi^D_{-2(k+m-1)+1-\frac{1}{2}}\dots \phi^D_{-1 -\frac{1}{2}}\phi^D_{-2(k-1)-\frac{1}{2}}\dots \phi^D_{2(m-1)-\frac{1}{2}}\phi^D_{2m-\frac{1}{2}}|0\rangle, \quad \text{for}\ k>-m.
\end{align*}
We  have
\begin{align*}
L^{\frac{1}{2}, \frac{m^2}{2}}_j \tilde{v}_{m, k}&=0, \quad \text{for \ any} \ j>0, k\geq -m;\\
L^{\frac{1}{2}, \frac{m^2}{2}}_0 \tilde{v}_{m, k}&=\frac{1}{4}(m+2k)^2\tilde{v}_{0, k}.
\end{align*}
\end{lem}
The proof is by direct calculation and we omit it.

\begin{prop}\label{propo:Virdecomp}
\textbf{Case I.} Let $\lambda =\frac{1}{2}$, $b$ is generic, i.e., $b\neq n_1 +\frac{m}{\sqrt{2}}$ for any $m, n_1\in \mathbb{Z}$. Then as Virasoro modules with central charge $c=1$
\begin{equation*}
\mathit{F_{\bar{0}}^{\ten \frac{1}{2}}}= \oplus_{\substack{n\in \mathbb{Z}\\ n\ \text{even}}}M\big(1, \frac{(b-n)^2}{2}\big);\quad
 \mathit{F_{\bar{1}}^{\ten \frac{1}{2}}}=\oplus_{\substack{n\in \mathbb{Z}\\ n\ \text{odd}}}M\big(1, \frac{(b-n)^2}{2}\big).
\end{equation*}
\textbf{Case II.}
 Let $\lambda =\frac{1}{2}$, $b=n_1 +\frac{m}{\sqrt{2}}$ for some unique $m, n_1\in \mathbb{Z}$ with $n_1$ \textbf{even}. Then as Virasoro modules with central charge $c=1$
\begin{align*}
\mathit{F_{\bar{0}}^{\ten \frac{1}{2}}} &=  \left(\oplus_{\substack{n\in \mathbb{Z}\\ n\ \text{even}\\ n\neq n_1}}M\big(1, \frac{(m-\sqrt{2}n)^2}{4}\big)\right)\bigoplus \left(\oplus_{\substack{n\geq 0\\ n\geq -m}}M\big(1, \frac{(m+2n)^2}{4}\big)\right);\\
\mathit{F_{\bar{1}}^{\ten \frac{1}{2}}}&=\oplus_{\substack{n\in \mathbb{Z}\\ n\ \text{odd}}}M\big(1, \frac{(m-\sqrt{2}n)^2}{4}\big).
\end{align*}
\textbf{Case III.}
 Let $\lambda =\frac{1}{2}$, $b=n_1 +\frac{m}{\sqrt{2}}$ for some unique $m, n_1\in \mathbb{Z}$ with $n_1$ \textbf{odd}. Then as Virasoro modules with central charge $c=1$
\begin{align*}
\mathit{F_{\bar{0}}^{\ten \frac{1}{2}}}&=\oplus_{\substack{n\in \mathbb{Z}\\ n\ \text{even}}}M\big(1, \frac{(m-\sqrt{2}n)^2}{4}\big);\\
\mathit{F_{\bar{1}}^{\ten \frac{1}{2}}} &=  \left(\oplus_{\substack{n\in \mathbb{Z}\\ n\  \text{odd}\\n\neq n_1}}M\big(1, \frac{(m-\sqrt{2}n)^2}{4}\big)\right)\bigoplus \left(\oplus_{\substack{n\geq 0\\ n\geq -m}}M\big(1, \frac{(m+2n)^2}{4}\big)\right).
\end{align*}
\end{prop}
We want to finish by showing an application of the above decomposition to calculating directly a positive sum (fermionic) representation of the characters of the Virasoro modules $M\left(\frac{1}{2}, \frac{1}{2}\right)$ and $M\left(\frac{1}{2}, 0\right)$. We can pick any $b$ and use  the Proposition above, together with the observation that we have a relation between the highest weights for $L^{1/2}_0$ and the highest weights for $L^{\lambda, b}_0 $, via the  connection represented by \eqref{eqn:h-normorder}.
 In particular for $\lambda =\frac{1}{2}$, $b=0$ we have
$L^{\frac{1}{2}, 0}(z^2)=L^{1}(z^2)$; and from  \eqref{eqn:h-normorder} we have that
\begin{equation}
L^{1/2}_0 =2 L^{1}_0 +\frac{1}{2}h^D_0.
\end{equation}
Hence we have from Proposition \ref{propo:Virdecomp}, Case II with $b=0$:
\begin{equation*}
ch_q M\big(\frac{1}{2}, \frac{1}{2}\big)= tr_{\mathit{F_{\bar{1}}^{\ten \frac{1}{2}}}} q^{L^{1/2}_0}= tr_{\mathit{F_{\bar{1}}^{\ten \frac{1}{2}}}} (q^2)^{L^{1}_0} q^{\frac{1}{2}h^D_0}=\sum_{\substack{n\in \mathbb{Z}\\ n\ \text{odd}}}ch_{q^2} M\big(1, \frac{n^2}{2}\big) q^{\frac{n}{2}}
=\sum_{\substack{n\in \mathbb{Z}\\ n\ \text{odd}}}\frac{q^{n^2+\frac{n}{2}}}{\prod_{i=1}^{\infty} (1-q^{2i})}.
\end{equation*}
Here we have used the well known character formula $ch_q M\big(1, \frac{n^2}{2}\big)=\frac{q^{n^2/2}}{\prod_{i=1}^{\infty} (1-q^{i})}$, for $n\in \mathbb{Z}, n\neq 0$, and we recover the formula \eqref{eqn:char1Vir} we obtained earlier.

Similarly, using Case II, $b=0$ and  the formula for $ch_q M\big(1, \frac{m^2}{4}\big)$ (see e.g. \cite{RochaCharidi}):
\[
ch_q M\big(1, \frac{m^2}{4}\big)=\frac{1}{\prod_{i=1}^{\infty} (1-q^{i})}\big(q^{m^2/4}-q^{(m+2)^2/4}\big), \quad \text{for} \ m\in \mathbb{Z},
\]
we get
\begin{align*}
ch_q M\big(\frac{1}{2}&, 0\big)= tr_{\mathit{F_{\bar{0}}^{\ten \frac{1}{2}}}} q^{L^{1/2}_0}= tr_{\mathit{F_{\bar{0}}^{\ten \frac{1}{2}}}} (q^2)^{L^{1}_0} q^{\frac{1}{2}h^D_0}=\sum_{\substack{n\in \mathbb{Z}\\ n\ \text{even}\\n\neq 0}}ch_{q^2} M\big(1, \frac{n^2}{2}\big) q^{\frac{n}{2}} +\sum_{\substack{n\geq 0}}ch_{q^2}M\big(1, n^2\big)\\
&=\sum_{\substack{n\in \mathbb{Z}\\ n\ \text{even}\\n\neq 0}}\frac{q^{n^2+\frac{n}{2}}}{\prod_{i=1}^{\infty} (1-q^{2i})}+ \sum_{\substack{n\geq 0}}\frac{q^{2n^2}-q^{2(n+1)^2}}{\prod_{i=1}^{\infty} (1-q^{2i})}=\sum_{\substack{n\in \mathbb{Z}\\ n\ \text{even}\\n\neq 0}}\frac{q^{n^2+\frac{n}{2}}}{\prod_{i=1}^{\infty} (1-q^{2i})}+ \frac{1}{\prod_{i=1}^{\infty} (1-q^{2i})}
\\
&=\sum_{\substack{n\in \mathbb{Z}\\ n\ \text{even}}}\frac{q^{n^2+\frac{n}{2}}}{\prod_{i=1}^{\infty} (1-q^{2i})}.
\end{align*}
Note that we could have used a generic $b$, for instance we have
\[
L^{\frac{1}{2}, -\frac{1}{4}}(z^2)=L^{1}(z^2) +\frac{1}{4z^3}h^D(z) +\frac{1}{32z^4}, \quad \quad L^{\frac{1}{2}, -\frac{1}{4}}_n =\frac{1}{2}L^{1/2}_{2n} +\frac{1}{32}\delta_{n, 0},
\]
but the resulting character formulas would have been the same.

We want to remark that this  application of multi-local  Virasoro field representations to calculating character formulas can be extended more generally to the discrete Virasoro series, not only to the Ising case of $c=\frac{1}{2}$. The main new idea is to  obtain a decomposition of (the vector space of) the irreducible modules represented by a  one-point local Virasoro field of discrete central charge $c=1-\frac{6}{(m+2)(m+3)}$ into irreducible modules represented by a multi-local Virasoro field of charge $1$ (or higher).

\section*{References}

\medskip


\def\cprime{$'$}

\end{document}